\documentclass[12pt]{article}
\usepackage{latexsym,amsmath,amssymb,amsfonts,amsthm,bbm}
\usepackage{natbib}
\usepackage{graphicx}
\usepackage{graphics}
\def\bx{\mathbf x}
\def\bz{\mathbf z}
\def\bs{\mathbf s}
\def\bX{\mathbf X}

\def\bR{\mathbb R}
\def\Fcal{\mathcal{F}}
\def\Pcal{\mathcal{P}}

\newtheorem{thm}{Theorem}
\newtheorem{prop}[thm]{Proposition}

\newtheorem{corollary}[thm]{Corollary}

\DeclareMathOperator*{\argmin}{argmin}
\DeclareMathOperator*{\argmax}{argmax}

\newenvironment{prooftitle}[1]{{\noindent \textsc{Proof #1}}\\}

\begin{document}

\title{Independent component analysis via nonparametric maximum likelihood estimation}
\author{Richard J. Samworth$^\ast$ and Ming Yuan$^\dag$\\
University of Cambridge and Georgia Institute of Technology}

\footnotetext[1]{
Statistical Laboratory, Centre for Mathematical Sciences, Wilberforce Road, Cambridge, CB3 0WB. The research of Richard J. Samworth was supported in part by a Leverhulme Research Fellowship and an EPSRC Early Career Fellowship.}

\footnotetext[2]{
H. Milton Stewart School of Industrial and Systems Engineering,
Georgia Institute of Technology, Atlanta, GA 30332. The research of Ming Yuan was supported in part by NSF Career Award DMS-0846234.}

\date{(\today)}

\maketitle
\begin{abstract}
Independent Component Analysis (ICA) models are very popular semiparametric models in which we observe independent copies of a random vector $X = AS$, where $A$ is a non-singular matrix and $S$ has independent components.  We propose a new way of estimating the unmixing matrix $W = A^{-1}$ and the marginal distributions of the components of $S$ using nonparametric maximum likelihood.  Specifically, we study the projection of the empirical distribution onto the subset of ICA distributions having log-concave marginals.  We show that, from the point of view of estimating the unmixing matrix, it makes no difference whether or not the log-concavity is correctly specified.  The approach is further justified by both theoretical results and a simulation study.
\end{abstract}

\newpage

\section{Introduction}
\label{sec:ontro}

In recent years, Independent Component Analysis (ICA) has seen an explosion in its popularity in diverse fields such as signal processing, machine learning, and medical imaging, to name a few. For a wide-ranging list of algorithms and applications of ICA, see the monograph by \citet{HKO2001}. In the ICA paradigm, one observes a random vector $X\in \bR^d$ that can be expressed as a non-singular linear transformation of $d$ mutually independent latent factors $S_1,\ldots,S_d$; thus $X=AS$ where $S=(S_1,\ldots, S_d)^{\sf T}$ and $A$ is a $d\times d$ full rank matrix often referred to as the mixing matrix. As such, ICA postulates the following model for the probability distribution $P$ of $X$: for any Borel set $B$ in $\bR^d$,
$$
P(B) = \prod_{j=1}^d P_j(w_j^{\sf T} B),
$$
where $W=(w_1,\ldots,w_d)^{\sf T}=A^{-1}$ is the so-called unmixing matrix, and $P_1,\ldots, P_d$ are the univariate probability distributions of the latent factors $S_1,\ldots, S_d$ respectively.

The goal of ICA, as in other blind source separation problems, is to infer from a sample $\bx_1,\ldots,\bx_n$ of independent observations of $X$, the independent factors $\bs_1=W\bx_1,\ldots,\bs_n=W\bx_n$, or equivalently the unmixing matrix $W$. This task is typically accomplished by first postulating a certain parametric family for the marginal probability distributions $P_1,\ldots, P_d$, and then optimising a contrast function involving $(W,P_1,\ldots,P_d)$. The contrast functions are often chosen to represent the mutual information as measured by Kullback--Leibler divergence or maximum entropy; or non-Gaussianity as measured by kurtosis or negentropy.  Alternatively, in recent years, methods for ICA have also been developed which assume $P_1,\ldots,P_d$ have smooth (log) densities, e.g. \citet{BachJordan2002}, \citet{HastieTibshirani2003}, \citet{SamTsy2004} and \citet{ChenBickel2006}.  Although more flexible than their aforementioned parametric peers, there remain unsettling questions about what happens if the smoothness assumptions on the marginal densities are violated, which may occur, in particular, when some of the marginal probability distributions $P_1,\ldots,P_d$ have atoms.  Another issue is that, in common with most other smoothing methods, a choice of tuning parameters is required to balance the fidelity to the observed data and the smoothness of the estimated marginal densities, and it is notoriously difficult to select these tuning parameters appropriately in practice.

In this paper, we argue that these assumptions and tuning parameters are unnecessary, and propose a new paradigm for ICA, based on the notion of nonparametric maximum likelihood, that is free of these burdens.  In fact, we show that the usual nonparametric (empirical) likelihood approach does not work in this context, and instead we proceed under the working assumption that the marginal distributions of $S_1,\ldots,S_d$ are log-concave. More specifically, we propose to estimate $W$ by maximising
$$
\log |\det W| + \frac{1}{n} \sum_{i=1}^n \sum_{j=1}^d \log f_j(w_j^{\sf T}\bx_i)
$$
over all $d\times d$ non-singular matrices $W=(w_1,\ldots,w_d)^{\sf T}$, and univariate log-concave densities $f_1,\ldots,f_d$. Remarkably, from the point of view of estimating the unmixing matrix $W$, it turns out that it makes no difference whether or not this hypothesis of log-concavity is correctly specified.  

The key to understanding how our approach works is to study what we call the log-concave ICA projection of a distribution on $\mathbb{R}^d$ onto the set of densities that satisfy the ICA model with log-concave marginals.  In Section~\ref{Sec:Notation} below, we define this projection carefully, and give necessary and sufficient conditions for it to make sense.  In Section~\ref{Sec:PdICA}, we prove that the log-concave projection of a distribution from the ICA model preserves both the ICA structure and the unmixing matrix.  Finally, in Section~\ref{Sec:Pd}, we derive a continuity property of log-concave ICA projections, which turns out to be important for understanding the theoretical properties of our ICA procedure. 

Our ICA estimating procedure uses the log-concave ICA projection of the empirical distribution of the data, and is studied in Section~\ref{Sec:EstProc}.  After explaining why the usual empirical likelihood approach cannot be used, we prove the consistency of our method.  We also present an iterative algorithm for the computation of our estimator. Our simulation studies in Section~\ref{Sec:Sim} confirm our theoretical results and show that the proposed method compares favourably with existing methods.

\section{Log-concave ICA projections}

Our proposed nonparametric maximum likelihood estimator can be viewed as the projection of the empirical distribution of $\bx_1,\ldots,\bx_n$ onto the space of ICA distributions with log-concave densities. To understand its behavior, it is useful to study the properties of such projections in general.

\subsection{Notation and overview}
\label{Sec:Notation}

Let $\mathcal{P}_k$ be the set of probability distributions $P$ on $\mathbb{R}^k$ satisfying $\int_{\mathbb{R}^k} \|x\| \, dP(x) < \infty$ and $P(H) < 1$ for all hyperplanes $H$, i.e. the probability measures in $\bR^k$ that have finite mean and are not supported in a translate of a lower dimensional linear subspace of $\bR^k$.  Here and throughout, $\|\cdot\|$ denotes the Euclidean norm on $\mathbb{R}^k$, and we will be interested in the cases $k=1$ and $k=d$.  Further, let $\mathcal{W}$ denote the set of non-singular $d \times d$ real matrices.  We use upper case letters to denote matrices in $\mathcal{W}$, and the corresponding lower case letters with subscripts to denote rows: thus $w_j^{\sf T}$ is the $j$th row of $W \in \mathcal{W}$.  Let $\mathcal{B}_k$ denote the class of Borel sets on $\mathbb{R}^k$. Then the ICA model $\mathcal{P}_d^{\mathrm{ICA}}$ is defined to be the set of $P \in \mathcal{P}_d$ of the form 
\begin{equation}
\label{Eq:PdICA}
P(B) = \prod_{j=1}^d P_j(w_j^{\sf T} B), \qquad \forall B \in \mathcal{B}_d,
\end{equation}
for some $W \in \mathcal{W}$ and $P_1,\ldots,P_d \in \mathcal{P}_1$.   As shown by \citet[Theorem~2.2]{DSS2011}, the condition $P \in \mathcal{P}_d$ is necessary and sufficient for the existence of a unique upper semi-continuous and log-concave density that is the closest to $P$ in the Kullback--Leibler sense. More precisely, let $\mathcal{F}_k$ denote the class of all upper semi-continuous, log-concave densities with respect to Lebesgue measure on $\mathbb{R}^k$. Then the projection $\psi^\ast: \Pcal_d \rightarrow \Fcal_d$ given by 
\[
\psi^*(P)=\argmax_{f \in \Fcal_d} \int_{\mathbb{R}^d} \log f \, dP
\]
is well-defined and surjective.  In what follows, we refer to $\psi^\ast$ as the log-concave projection operator and $f^\ast:=\psi^\ast(P)$ as the log-concave projection of $P$.  By a slight abuse of notation, we also use $\psi^*$ to denote the log-concave projection from $\Pcal_1$ to $\Fcal_1$. 

Although the log-concave projection operator does play a role in this paper, our main interest is in a different projection, onto the subset of $\Fcal_d$ consisting of those densities that satisfy the ICA model.  This class is given by
\begin{equation}
\label{Eq:FdICA}
\mathcal{F}_d^{\mathrm{ICA}} = \biggl\{f \in \mathcal{F}_d: f(x) = |\det W|\prod_{j=1}^d f_j(w_j^{\sf T} x) \, dx \text{ for some } W \in \mathcal{W} \text{ and } f_1,\ldots,f_d \in \mathcal{F}_1\biggr\}.
\end{equation}
Note that, in this representation, if $X$ has density $f \in \mathcal{F}_d^{\mathrm{ICA}}$, then $w_j^{\sf T} X$ has density $f_j$.  The corresponding log-concave ICA projection operator $\psi^{**}(\cdot)$ is defined for any distribution $P$ on $\mathbb{R}^d$ by
\[
\psi^{**}(P)=\argmax_{f \in \Fcal_d^{\mathrm{ICA}}} \int_{\mathbb{R}^d} \log f \, dP.
\]
We also write $L^{**}(P) = \sup_{f \in \Fcal_d^{\mathrm{ICA}}} \int_{\mathbb{R}^d} \log f \, dP$.  
\begin{prop}
\label{Prop:Cases}
\begin{enumerate}
\item If $\int_{\mathbb{R}^d} \|x\| \, dP(x) = \infty$, then $L^{**}(P) = -\infty$ and $\psi^{**}(P) = \Fcal_d^{\mathrm{ICA}}$. 
\item If $\int_{\mathbb{R}^d} \|x\| \, dP(x) < \infty$, but $P(H) = 1$ for some hyperplane $H$, then $L^{**}(P) = \infty$ and $\psi^{**}(P) = \emptyset$.
\item If $P \in \mathcal{P}_d$, then $L^{**}(P) \in \mathbb{R}$ and $\psi^{**}(P)$ defines a non-empty, proper subset of $\Fcal_d^{\mathrm{ICA}}$.
\end{enumerate}
\end{prop}
In view of Proposition~\ref{Prop:Cases}, and to avoid lengthy discussion of trivial exceptional cases, we henceforth consider $\psi^{**}(\cdot)$ as being defined on $\mathcal{P}_d$.  In contrast to $\psi^*(P)$, which defines a unique element of $\Fcal_d$, the log-concave ICA projection operator $\psi^{**}(P)$ may not define a unique element of $\Fcal_d^{\mathrm{ICA}}$, even for $P \in \mathcal{P}_d$.  For instance, consider the situation where $P$ is the uniform distribution on the closed unit disk in $\mathbb{R}^2$ equipped with the Euclidean norm.  Here, the spherical symmetry means that the choice of $W \in \mathcal{W}$ is arbitrary.  In fact, after a straightforward calculation, it can be shown that $\psi^{**}(P)$ consists of those $f \in \mathcal{F}_d^{\mathrm{ICA}}$ where, in the representation~(\ref{Eq:FdICA}), $W \in \mathcal{W}$ is arbitrary and $f_1,f_2 \in \mathcal{F}_1$ are given by $f_1(x) = f_2(x) = \frac{2}{\pi}(1-x^2)^{1/2}\mathbbm{1}_{\{x \in [-1,1]\}}$.  It is certainly possible to make different choices of $W$ that yield different elements of $\mathcal{F}_d^{\mathrm{ICA}}$.  This example shows that, in general, we must think of $\psi^{**}(P)$ as defining a subset of $\mathcal{F}_d^{\mathrm{ICA}}$.

The relationship between the spaces introduced above and the projection operators is illustrated in the diagram below:
\[
\begin{array}{lll}
\mathcal{P}_d & \stackrel{\psi^*}{\longrightarrow} & \mathcal{F}_d \\
                           & \stackrel{\psi^{**}}{\searrow} & \\
 \mathcal{P}_d^{\mathrm{ICA}} & \stackrel{\psi^{**}|_{\mathcal{P}_d^{\mathrm{ICA}}}}{\longrightarrow} &  \mathcal{F}_d^{\mathrm{ICA}}
\end{array}
\]
Our next subsection studies the restriction of $\psi^{**}$ to $\mathcal{P}_d^{\mathrm{ICA}}$, denoted $\psi^{**}|_{\mathcal{P}_d^{\mathrm{ICA}}}$; Section~\ref{Sec:PdICA} examines $\psi^{**}$ more generally as a map on $\mathcal{P}_d$.

\subsection{Log-concave projections of the ICA model}
\label{Sec:PdICA}

Our first result in this subsection characterises $\psi^{**}|_{\mathcal{P}_d^{\mathrm{ICA}}}$.
\begin{thm}
\label{Thm:PdICA}
If $P \in \mathcal{P}_d^{\mathrm{ICA}}$, then $\psi^{**}(P)$ defines a unique element of $\Fcal_d^{\mathrm{ICA}}$.  The map $\psi^{**}|_{\mathcal{P}_d^{\mathrm{ICA}}}$ is surjective, and coincides with $\psi^*|_{\mathcal{P}_d^{\mathrm{ICA}}}$.  Moreover, suppose that $P \in \mathcal{P}_d^{\mathrm{ICA}}$, so that
\[
P(B) = \prod_{j=1}^d P_j(w_j^{\sf T} B), \quad \forall B \in \mathcal{B}_d,
\]
for some $W \in \mathcal{W}$ and $P_1,\ldots,P_d \in \mathcal{P}_1$.  Then $f^{**} = \psi^{**}(P)$ can be written as 
\[
f^{**}(x) = |\det W| \prod_{j=1}^d f_j^*(w_j^{\sf T} x),
\]
where $f_j^* = \psi^*(P_j)$.
\end{thm}
It is interesting to observe that log-concave projection operator $\psi^*$ preserves the ICA structure.  But perhaps the most important aspect of this result is the fact that the same unmixing matrix $W$ can be used to represent both the original ICA model and its log-concave projection.  This observation lies at the heart of the rationale for our approach to ICA.  

A remaining concern is that the unmixing matrix may not be identifiable.  For instance, applying the same permutation to the rows of $W$ and the vector of marginal distributions $(P_1,\ldots,P_d)$ leaves the distribution unchanged; similarly, the same effect occurs if we multiply any of the rows of $W$ by a scaling factor and applying the corresponding scaling factor to the relevant marginal distribution.  The question of identifiability for ICA models was first addressed by \citet{Comon1994}, who assumed that $W$ is orthogonal, and was settled in the general case by \citet{ErikssonKoivunen2004}.  One way to state their result is as follows: suppose that a probability measure $P$ on $\mathbb{R}^d$ has two representations as
\begin{equation}
\label{Eq:TwoReps}
P(B) = \prod_{j=1}^d P_j(w_j^{\sf T} B) = \prod_{j=1}^d \tilde{P}_j(\tilde{w}_j^{\sf T} B) \quad \forall B \in \mathcal{B}_d,
\end{equation}
where $W$, $\tilde{W} \in \mathcal{W}$ and $P_1,\ldots,P_d,\tilde{P}_1,\ldots,\tilde{P}_d$ are probability measures on $\mathbb{R}$.  Then the pair of conditions that $P_1,\ldots,P_d$ are not Dirac point masses and not more than one of $P_1,\ldots,P_d$ is Gaussian is necessary and sufficient for the existence of a permutation $\pi$ of $\{1,\ldots,d\}$ and scaling vector $\epsilon = (\epsilon_1,\ldots,\epsilon_d) \in (\mathbb{R} \setminus \{0\})^d$ such that $\tilde{P}_j(B_j) = P_{\pi(j)}(\epsilon_j B_j)$ for all $B_j \in \mathcal{B}_1$, and $\tilde{w}_j = \epsilon_j^{-1} w_{\pi(j)}$.  When such a permutation and scaling factor exist for any two ICA representations of $P$, we say that \emph{the ICA representation of $P$ is identifiable}, or simply that \emph{P is identifiable}.  By analogy, we define $f \in \mathcal{F}_d^{\mathrm{ICA}}$ to be identifiable if not more than one of $f_1,\ldots,f_d$ in the representation~(\ref{Eq:FdICA}) is Gaussian.

Our next result shows that $\psi^{**}$ preserves the identifiability of the ICA model.  Together with Theorem~\ref{Thm:PdICA}, we see that if $P \in \mathcal{P}_d^{\mathrm{ICA}}$ is identifiable, then the unmixing matrices of $P \in \mathcal{P}_d^{\mathrm{ICA}}$ and $\psi^{**}(P)$ are identical up to the permutation and scaling transformations described above.
\begin{thm}
\label{Thm:Identifiability}
Let $P \in \mathcal{P}_d^{ICA}$.  Then $\psi^{**}(P)$ is identifiable if and only if $P$ is identifiable. 
\end{thm}

\subsection{General log-concave ICA projections}
\label{Sec:Pd}

We now consider the general log-concave ICA projection $\psi^{**}$ defined on $\mathcal{P}_d$.  Define the Mallows distance $d$ (also known as the Wasserstein distance) between probability measures $P$ and $Q$ on $\mathbb{R}^d$ with finite mean by
\[
d(P,\tilde{P}) = \inf_{(X,\tilde{X}) \sim (P,\tilde{P})} \mathbb{E}\|X-Y\|,
\]
where the infimum is taken over all pairs $(X,Y)$ of random vectors $X \sim P$ and $\tilde{X} \sim \tilde{P}$ on a common probability space.  Recall that $d(P^n,P) \rightarrow 0$ if and only if both $P^n \stackrel{d}{\rightarrow} P$ and $\int_{\mathbb{R}^d} \|x\| \, dP^n(x) \rightarrow \int_{\mathbb{R}^d} \|x\| \, dP(x)$.  We are interested in the continuity of $\psi^{**}$.  
\begin{prop}
\label{Prop:Pd}
Let $P, P^1,P^2,\ldots$ be probability measures in $\mathcal{P}_d$ with $d(P^n,P) \rightarrow 0$ as $n \rightarrow \infty$.  Then $L^{**}(P^n) \rightarrow L^{**}(P)$.  Moreover, 
\[
\sup_{f^n \in \psi^{**}(P^n)} \inf_{f \in \psi^{**}(P)} \int_{\mathbb{R}^d} |f^n-f| \rightarrow 0
\]
as $n \rightarrow \infty$.
\end{prop}
The second part of this proposition says that any element of $\psi^{**}(P^n)$ is arbitrarily close in total variation distance to some element of $\psi^{**}(P)$ once $n$ is sufficiently large.  In the special case where $\psi^{**}(P)$ consists of only a single element, we can say more.  It is convenient to let $\Pi_d$ denote the set of permutations of $\{1,\ldots,d\}$, and write $(W,f_1,\ldots,f_d) \stackrel{\mathrm{ICA}}{\sim} f$ if $W \in \mathcal{W}$ and $f_1,\ldots,f_d \in \mathcal{F}_1$ can be used to give an ICA representation of $f \in \mathcal{F}_d^{\mathrm{ICA}}$ in (\ref{Eq:FdICA}).  Similarly, we write $(W,P_1,\ldots,P_d) \stackrel{\mathrm{ICA}}{\sim} P$ if $W \in \mathcal{W}$ and $P_1,\ldots,P_d \in \mathcal{P}_1$ represent $P \in \mathcal{P}_d^{\mathrm{ICA}}$ in (\ref{Eq:PdICA}).
\begin{thm}
\label{Thm:Conv}
Suppose that $P \in \mathcal{P}_d^{\mathrm{ICA}}$, and write $f^{**} = \psi^{**}(P)$.  If $P^1,P^2,\ldots \in \mathcal{P}_d$ are such that $d(P^n,P) \rightarrow 0$, then
\[
\sup_{f^n \in \psi^{**}(P^n)} \int_{\mathbb{R}^d} |f^n - f^{**}| \rightarrow 0.
\]
Suppose further that $P$ is identifiable and that $(W,P_1,\ldots,P_d) \stackrel{\mathrm{ICA}}{\sim} P$.  Then
\begin{align*}
\sup_{f^n \in \psi^{**}(P^n)} \sup_{(W^n,f_1^n,\ldots,f_d^n) \stackrel{\mathrm{ICA}}{\sim} f^n} \inf_{\pi^n \in \Pi_d} \inf_{\epsilon_1^n,\ldots,\epsilon_d^n \in \mathbb{R} \setminus \{0\}} \biggl\{\|&(\epsilon_j^n)^{-1} w_{\pi^n(j)}^n - w_j\| \\
&+ \int_{-\infty}^\infty \bigl| |\epsilon_j^n| f_{\pi^n(j)}^n(\epsilon_j^n x) - f_j^*(x)\bigr| \, dx\biggr\} \rightarrow 0,
\end{align*}
for each $j=1,\ldots,d$, where $f_j^* = \psi^*(P_j)$.  As a consequence, for sufficiently large $n$, every $f^n \in \psi^{**}(P^n)$ is identifiable.
\end{thm}

The first part of Theorem~\ref{Thm:Conv} show that if $P \in P_d^{\mathrm{ICA}}$ and $\tilde{P} \in P_d$ are close in Mallows distance, then every $\tilde{f} \in \psi^{**}(\tilde{P})$ is close to the corresponding (unique) log-concave ICA projection $f = \psi^{**}(P)$ in total variation distance.  The second part shows further that if $P$ is identifiable, then up to permutation and scaling, every $\tilde{f} \in \psi^{**}(\tilde{P})$ and every choice of unmixing matrix $\tilde{W}$ and marginal densities $\tilde{f}_1,\ldots,\tilde{f}_d$ in the ICA representation of $\tilde{f}$ is close to the unmixing matrix $W$ and marginal densities $f_1,\ldots,f_d$ in the ICA representation of $f$. 

To conclude this subsection, we remark that, by analogy with the situation when $P \in \mathcal{P}_d^{\mathrm{ICA}}$ described in Theorem~\ref{Thm:PdICA}, if $P \in \mathcal{P}_d$ and $X \sim P$, any $f^{**} \in \psi^{**}(P)$ can be written as
\[
f^{**}(x) = |\det W| \prod_{j=1}^d f_j^*(w_j^{\sf T} x),
\]
for some $W \in \mathcal{W}$, where $f_j^* = \psi^*(P_j)$, and $P_j$ is the marginal distribution of $w_j^{\sf T} X$.  This observation reduces the maximisation problem involved in computing $\psi^{**}(P)$ to a finite-dimensional one (over $W \in \mathcal{W}$), and follows because
\begin{align*}
\sup_{f \in \Fcal_d^{\mathrm{ICA}}} \int_{\mathbb{R}^d} \log f \, dP &= \sup_{W \in \mathcal{W}} \sup_{f_1,\ldots,f_d \in \Fcal_1} \biggl\{\log |\det W| + \sum_{j=1}^d \int_{\mathbb{R}^d} \log f_j(w_j^{\sf T} x) \, dP(x)\biggr\} \\
&= \sup_{W \in \mathcal{W}} \biggl\{\log |\det W| + \sum_{j=1}^d \int_{\mathbb{R}^d} \log f_j^*(w_j^{\sf T} x) \, dP(x)\biggr\}.
\end{align*}

\section{Nonparametric maximum likelihood estimation for ICA models}
\label{Sec:EstProc}

We are now in position to study the proposed nonparametric maximum likelihood estimator.
\subsection{Estimating procedure and theoretical properties}

Now assume $\bx_1,\bx_2,\ldots$ are independent copies of a random vector $X \in \mathbb{R}^d$ satisfying the ICA model.  Thus $X = AS$, where $A = W^{-1} \in \mathcal{W}$ and $S = (S_1,\ldots,S_d)^{\sf T}$ has independent components.  In this section, we study a nonparametric maximum likelihood estimator of $W$ and the marginal distributions $P_1,\ldots,P_d$ of $S_1,\ldots,S_d$ based on $\bx_1,\ldots,\bx_n$, where $n \geq d+1$.

We start by noting that the usual nonparametric maximum likelihood estimate does not work.  Indeed, in the spirit of empirical likelihood \citep{Owen1990}, it would suffice to consider, for a given $W = (w_1,\ldots,w_d)^{\sf T} \in \mathcal{W}$, estimates $\tilde{P}_j$ of the marginal distribution $P_j$, supported on $w_j^{\sf T}\bx_1,\ldots,w_j^{\sf T}\bx_n$.  This leads to the nonparametric likelihood
\begin{equation}
\label{Eq:NonLike}
L(W,\tilde{P}_1,\ldots,\tilde{P}_d) = \prod_{i=1}^n \prod_{j=1}^d \tilde{p}_{ij},
\end{equation}
where $\tilde{p}_{ij} = \tilde{P}_j(w_j^{\sf T}\bx_i)$.  Let $J$ denote a subset of $(d+1)$ distinct indices in $\{1,\ldots,n\}$, and let $\bX_J$ denote the $d \times (d+1)$ matrix obtained by extracting the columns of $\bX = (\bx_1,\ldots,\bx_n)$ with indices in $J$.   Now let $\bX_{(-j)}$ denote the $d \times d$ matrix obtained by removing the $j$th column of $\bX_J$.  Let $W_J \in \mathcal{W}$ have $j$th row $w_j = (\bX_{(-j)}^{-1})^{\sf T}\mathbf{1}_d$, for $j=1,\ldots,d$, where $\mathbf{1}_d$ is a $d$-vector of ones.  Our next result shows that every $W_J$ corresponds to a maximiser of the nonparametric likelihood~(\ref{Eq:NonLike}).   
\begin{prop}
\label{Prop:EmpLike}
Suppose that $\bx_1,\ldots,\bx_n$ are in general position.  Then for any choice $J$ of $(d+1)$ distinct indices in $\{1,\ldots,n\}$, there exist $\hat{P}_1,\ldots,\hat{P}_d \in \mathcal{P}_1$ such that $(W_J, \hat{P}_1,\ldots,\hat{P}_d)$ maximises $L(\cdot)$.
\end{prop}
If $X$ has a density with respect to Lebesgue measure on $\mathbb{R}^d$, then with probability 1, every subset of $\bx_1,\ldots,\bx_n$ of size $(d+1)$ is in general position.  On the other hand, there is no reason for different choices of $J$ to yield similar estimates $W_J$, so we cannot hope for such an empirical likelihood-based procedure to be consistent.

As a remedy, we propose to estimate $P^0 \in \mathcal{P}_d^{ICA}$ by $\psi^{**}(\hat{P}^n)$, where $\hat{P}^n$ denotes the empirical distribution of $\bx_1,\ldots,\bx_n \sim P^0$.  More explicitly, we estimate the unmixing matrix and the marginals by maximising the log-likelihood
\begin{equation}
\label{Eq:LogLike}
\ell^n(W,f_1,\ldots,f_d) = \ell^n(W,f_1,\ldots,f_d;\bx_1,\ldots,\bx_n) = \log |\det W| + \frac{1}{n} \sum_{i=1}^n \sum_{j=1}^d \log f_j(w_j^{\sf T}\bx_i)
\end{equation}
over $W \in \mathcal{W}$ and $f_1,\ldots,f_d \in \mathcal{F}_1$.  Note from Proposition~\ref{Prop:Cases} that $\psi^{**}(\hat{P}^n)$ exists as a proper subset of $\mathcal{F}_d^{\mathrm{ICA}}$ once the convex hull of $\bx_1,\ldots,\bx_n$ is $d$-dimensional, which happens with probability 1 for sufficiently large $n$.  As a direct consequence of Theorem~\ref{Thm:Conv} and the fact that $d(\hat{P}^n,P^0) \stackrel{a.s.}{\rightarrow} 0$, we have the following consistency result.
\begin{corollary}
Suppose that $P^0 \in \mathcal{P}_d^{ICA}$ is identifiable and is represented by $W^0 \in \mathcal{W}$ and $P_1^0,\ldots,P_d^0 \in \mathcal{P}_1$.  Then for any maximiser $(\hat{W}^n,\hat{f}_1^n,\ldots,\hat{f}_d^n)$ of $\ell^n(W,f_1,\ldots,f_d)$ over $W \in \mathcal{W}$ and $f_1,\ldots,f_d \in \Fcal_1$, there exist a permutation $\hat{\pi}^n$ of $\{1,\ldots,d\}$ and scaling factors $\hat{\epsilon}_1^n,\ldots,\hat{\epsilon}_d^n \in \mathbb{R}\setminus  \{0\}$ such that
\[
(\hat{\epsilon}_j^n)^{-1}\hat{w}_{\hat{\pi}^n(j)}^n \stackrel{a.s.}{\rightarrow} w_j^0 \quad \text{and} \quad \int_{-\infty}^\infty \bigl| |\hat{\epsilon}_j^n|\hat{f}_{\hat{\pi}^n(j)}^n(\hat{\epsilon}_j^n x) - f_j^*(x)\bigr| \, dx \stackrel{a.s.}{\rightarrow} 0,
\]
for $j=1,\ldots,d$, where $f_j^* = \psi^*(P_j^0)$. 
\end{corollary}

\subsection{Pre-whitening}
\label{Sec:Pre-Whiten}

Pre-whitening is a standard pre-processing technique in the ICA literature; see \citet[][pp.140--141]{HKO2001} or \citet{ChenBickel2005}.  In this subsection, we explain the rationale for pre-whitening and the simplifications it provides.

Assume for now that $P \in \mathcal{P}_d^{\mathrm{ICA}}$ and $\int_{\mathbb{R}^d} \|x\|^2 \, dP(x) < \infty$, and let $\Sigma$ denote the (positive-definite) covariance matrix corresponding to $P$.  Consider the ICA model $X = AS$, where $X \sim P$, the mixing matrix $A$ is non-singular and $S = (S_1,\ldots,S_d)$ has independent components with $S_j \sim P_j$.  Assuming without loss of generality that each component of $S$ has unit variance, we can write $\Sigma^{-1/2}X = \Sigma^{-1/2}AS \equiv \tilde{A}S$, say, where $\tilde{A}$ belongs to the set $O(d)$ of orthogonal $d \times d$ matrices.  Thus the unmixing matrix $W$ belongs to the set $O(d)\Sigma^{-1/2} = \{O\Sigma^{-1/2}:O \in O(d)\}$.

It follows that, if $\Sigma$ were known, we could maximise $\ell^n$ with the restriction that $W \in O(d)\Sigma^{-1/2}$.  In practice, $\Sigma$ is typically unknown, but we can estimate it using the sample covariance matrix $\hat{\Sigma}$.  For $n$ large enough that the convex hull of $\bx_1,\ldots,\bx_n$ is $d$-dimensional, we can therefore consider maximising
\[
\ell^n(W,f_1,\ldots,f_d;\bx_1,\ldots,\bx_n)
\]
over $W \in O(d)\hat{\Sigma}^{-1/2}$ and $f_1,\ldots,f_d \in \mathcal{F}_1$.  Denote such a maximiser by $(\hat{\hat{W}}^n,\hat{\hat{f}}_1^n,\ldots,\hat{\hat{f}}_d^n)$.  The corollary below shows that, under a second moment condition, $\hat{\hat{W}}^n$ and $\hat{\hat{f}}_1^n,\ldots,\hat{\hat{f}}_d^n$ have the same asymptotic properties as the original estimators $\hat{W}^n$ and $\hat{f}_1^n,\ldots,\hat{f}_d^n$.
\begin{corollary}
\label{Cor:Hats}
Suppose that $P^0 \in \mathcal{P}_d^{\mathrm{ICA}}$ is identifiable, is represented by $W^0 \in \mathcal{W}$ and $P_1^0,\ldots,P_d^0 \in \mathcal{P}_1$ and that $\int_{\mathbb{R}^d} \|x\|^2 \, dP^0(x) < \infty$.  Then with probability 1 for sufficiently large $n$, a maximiser $(\hat{\hat{W}}^n,\hat{\hat{f}}_1^n,\ldots,\hat{\hat{f}}_d^n)$ of $\ell^n(W,f_1,\ldots,f_d)$ over $W \in O(d)\hat{\Sigma}^{-1/2}$ and $f_1,\ldots,f_d \in \Fcal_1$ exists.  Moreover, for any such maximiser, there exist a permutation $\hat{\hat{\pi}}^n$ of $\{1,\ldots,d\}$ and scaling factors $\hat{\hat{\epsilon}}_1^n,\ldots,\hat{\hat{\epsilon}}_d^n \in \mathbb{R}\setminus  \{0\}$ such that
\[
(\hat{\hat{\epsilon}}_j^n)^{-1}\hat{\hat{w}}_{\hat{\hat{\pi}}^n(j)}^n \stackrel{a.s.}{\rightarrow} w_j^0 \quad \text{and} \quad \int_{-\infty}^\infty \bigl| |\hat{\hat{\epsilon}}_j^n|\hat{\hat{f}}_{\hat{\hat{\pi}}^n(j)}^n(\hat{\hat{\epsilon}}_j^n x) - f_j^*(x)\bigr| \, dx \stackrel{a.s.}{\rightarrow} 0,
\]
where $f_j^* = \psi^*(P_j^0)$.
\end{corollary}
An alternative, equivalent way of computing $(\hat{\hat{W}}^n,\hat{\hat{f}}_1^n,\ldots,\hat{\hat{f}}_d^n)$ is to \emph{pre-whiten} the data by replacing $\bx_1,\ldots,\bx_n$ with $\bz_1 = \hat{\Sigma}^{-1/2}\bx_1,\ldots,\bz_n = \hat{\Sigma}^{-1/2}\bx_n$, and then maximise 
\[
\ell^n(O,g_1,\ldots,g_d;\bz_1,\ldots,\bz_n)
\]
over $O \in O(d)$ and $g_1,\ldots,g_d \in \Fcal_1$.  If $(\hat{O}^n,\hat{g}_1^n,\ldots,\hat{g}_d^n)$ is such a maximiser, we can then set $\hat{\hat{W}}^n = \hat{O}^n\hat{\Sigma}^{-1/2}$ and $\hat{\hat{f}}_j^n = \hat{g}_j^n$.  Note that pre-whitening breaks down the estimation of the $d^2$ parameters in $W$ into two stages: first, we use $\hat{\Sigma}$ to estimate the $d(d+1)/2$ free parameters of the symmetric, positive definite matrix $\Sigma$, leaving only the maximisation over the $d(d-1)/2$ free parameters of $O \in O(d)$ at the second stage.  The advantage of this approach is that it facilitates more stable maximisation algorithms, such as the one described in the next subsection.

\subsection{Computational algorithm}

In this subsection, we address the challenge of maximising 
\[
\ell^n(W,g_1,\ldots,g_d;\bz_1,\ldots,\bz_n)
\]
over $W \in O(d)$ and $g_1,\ldots,g_d \in \Fcal_1$.   As a starting point, we choose $W$ to be randomly distributed according to Haar measure on the set $O(d)$ of $d \times d$ orthogonal matrices.  A simple way of generating $W$ with this distribution is to generate a $d \times d$ matrix $Z$ whose entries are independent $N(0,1)$ random variables, compute the $QR$-factorisation $Z = QR$, and let $W = Q$.

Our proposed algorithm then alternates between maximising the log-likelihood over $f_1,\ldots,f_d$ for fixed $W$, and then over $W$ for fixed $f_1,\ldots,f_d$.  The first of these steps is straightforward given Theorem~\ref{Thm:PdICA} and the recent work on log-concave density estimation: we set $f_j$ to be the log-concave maximum likelihood estimator of the data $w_j^{\sf T}\bx_1,\ldots,w_j^{\sf T}\bx_n$.  This can be computed using the \texttt{R} package \texttt{logcondens} \citep{RufibachDuembgen2006,DuembgenRufibach2011}. 

This leaves the challenge of updating $W \in O(d)$.  In order to describe our proposal, we recall some basic facts from differential geometry.  The set $O(d)$ is a $d(d-1)/2$-dimensional submanifold of $\mathbb{R}^{d^2}$.  The tangent space at $W \in O(d)$ is $T_W O(d) := \{WY:Y = -Y^{\sf T}\}$.  In fact, if we define the natural inner product $\langle \cdot , \cdot \rangle$ on $T_W O(d) \times T_W O(d)$ by $\langle U,V\rangle = \mathrm{tr} (UV^{\sf T})$, then $O(d)$ becomes a Riemannian manifold.  (Note that if we think of $U$ and $V$ as vectors in $\mathbb{R}^{d^2}$, then this inner product is simply the Euclidean inner product.)  

There is no loss of generality in assuming $W$ belongs to the Riemannian manifold $SO(d)$, the set of special orthogonal matrices having determinant 1.  We can now define geodesics on $SO(d)$, recalling that the matrix exponential is given by
\[
\exp(Y) = I + \sum_{r=1}^\infty \frac{Y^r}{r!}.
\]
The unique geodesic passing through $W \in SO(d)$ with tangent vector $WY$ (where $Y = -Y^{\sf T}$) is the map $\alpha:[0,1] \rightarrow SO(d)$ given by $\alpha(t) = W\exp(tY)$.  

We update $W$ by moving along a geodesic in $SO(d)$, but need to choose an appropriate skew-symmetric matrix $Y$, which ideally should (at least locally) give a large increase in the log-likelihood.  The key to finding such a direction is Proposition~\ref{Prop:Diff} below.  To set the scene for this result, observe that for $x \in [\min(w_j^{\sf T}\bx_1,\ldots,w_j^{\sf T}\bx_n),\max(w_j^{\sf T}\bx_1,\ldots,w_j^{\sf T}\bx_n)]$, we can write 
\begin{equation}
\label{Eq:logfj}
\log f_j(x) = \min_{k=1,\ldots,m_j} (b_{jk}x - \beta_{jk}),
\end{equation}
for some $b_{jk},\beta_{jk} \in \mathbb{R}$ \citep[e.g.][]{CSS2010}.  Since we may assume that $b_{j1},\ldots,b_{jm_j}$ are strictly decreasing, the minimum in~\eqref{Eq:logfj} is attained in either one or two indices.  It is convenient to let $\mathcal{K}_{ij} = \argmin_{k=1,\ldots,m_j} (b_{jk}w_j^{\sf T}\bx_i - \beta_{jk})$.
\begin{prop}
\label{Prop:Diff}
Consider the map $g:SO(d) \rightarrow \mathbb{R}$ given by 
\[
g(W) = \frac{1}{n} \sum_{i=1}^n \sum_{j=1}^d \min_{k=1,\ldots,m_j} (b_{jk}w_j^{\sf T}\bx_i - \beta_{jk}).
\]
Let $Y$ be a skew-symmetric matrix and let $c_j$ denote the $j$th row of $WY$.  If $|\mathcal{K}_{ij}| = 1$, let $k_{ij}$ denote the unique element of $\mathcal{K}_{ij}$.  If $|\mathcal{K}_{ij}| = 2$, write $\mathcal{K}_{ij} = \{k_{ij1},k_{ij2}\}$.  If $c_j^{\sf T}\bx_i \geq 0$, let $k_{ij} = k_{ijl}$, where $l = \argmin_{l = 1,2} b_{k_{ijl}}$; if $c_j^{\sf T}\bx_i < 0$, let $k_{ij} = k_{ijl}$, where $l = \argmax_{l=1,2} b_{k_{ijl}}$.  Then the one-sided directional derivative of $g$ at $W$ in the direction $WY$ is 
\[
\nabla_{WY}g(W) := \frac{1}{n} \sum_{i=1}^n \sum_{j=1}^d b_{jk_{ij}} c_j^{\sf T} \bx_i.
\]
\end{prop}
For $1 < s < r < d$, let $Y_{r,s}$ denote the $d \times d$ matrix with $Y_{r,s}(r,s) = 1/\sqrt{2}$, $Y_{r,s}(s,r) = -1/\sqrt{2}$ and all other entries equal to zero.  Then $\mathcal{Y}^+ = \{Y_{r,s}: 1 < s < r < d\}$ forms an orthonormal basis for the set of skew-symmetric matrices.  Let $\mathcal{Y}^- = \{-Y:Y \in \mathcal{Y}^+\}$.  We choose $Y^{\max} \in \mathcal{Y}^+ \cup \mathcal{Y}^-$ to maximise $\nabla_{WY}g(Y)$.

We therefore update $W$ with $W \exp(\epsilon Y^{\max})$, and it remains to select $\epsilon$.  This we propose to choose by means of a backtracking line search.   Specifically, we fix $\alpha \in (0,1)$ and $\epsilon = 1$, and if
\begin{equation}
\label{Eq:LineSearch}
g(W \exp(\epsilon Y^{\max})) > g(W) + \alpha \epsilon \nabla_{WY^{\max}}g(W),
\end{equation}
we accept a move from $W$ to $W \exp(\epsilon Y^{\max})$.  Otherwise, we successively reduce $\epsilon$ by a factor of $\gamma \in (0,1)$ until \eqref{Eq:LineSearch} is satisfied, and then move to $W \exp(\epsilon Y^{\max})$.  In our implementation, we used $\alpha = 0.3$ and $\gamma = 1/2$.

Our algorithm produces a sequence $(W^{(1)},f_1^{(1)},\ldots,f_d^{(1)}),(W^{(2)},f_1^{(2)},\ldots,f_d^{(2)}),\ldots$.  We terminate the algorithm once
\[
\frac{\ell^n(W^{(t)},f_1^{(t)},\ldots,f_d^{(t)}) - \ell^n(W^{(t-1)},f_1^{(t-1)},\ldots,f_d^{(t-1)})}{|\ell^n(W^{(t-1)},f_1^{(t-1)},\ldots,f_d^{(t-1)})|} < \eta,
\]
where, in our implementation, we chose $\eta = 10^{-7}$.  As with other ICA algorithms, global convergence is not guaranteed, so we used 10 random starting points and took the solution with the highest log-likelihood. 

\section{Numerical Experiments}
\label{Sec:Sim}

To illustrate the practical merits of our proposed nonparametric maximum likelihood estimation method for ICA models, we conducted several sets of numerical experiments.  To fix ideas, we focus on two-dimensional signals, that is $d=2$.  The components of the signal were generated independently, and then rotated by $\pi/3$, so the mixing matrix is
\[
A = \begin{pmatrix} 1/2 & -\sqrt{3}/2 \\ \sqrt{3}/2 & 1/2\end{pmatrix}.
\]
Our goal is to reconstruct the signal and estimate $A$, or equivalently $W=A^{-1}$, based on $n=200$ observations of the rotated input.

\begin{figure}[htbp]
\begin{center}
\includegraphics[height=0.7\textheight,width=0.5\textwidth,angle=270]{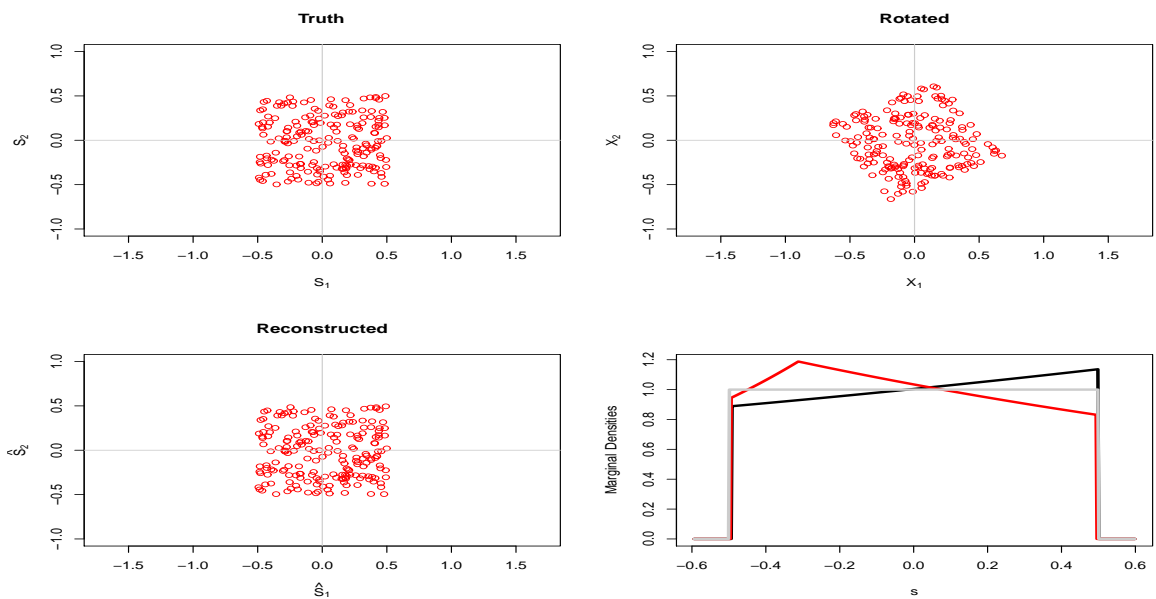}
\caption{Uniform signal: Top left panel, top right panel and bottom left panel give the true signal, rotated observations and the reconstructed signal respectively. The bottom right panel gives the estimated marginal densities along with the true marginal (grey line).}
\label{fig:uniform}
\end{center}
\end{figure}

We first consider a typical example in the ICA literature where the density of each component of the true signal is uniform on the interval $[-0.5,0.5]$.  The top left panel of Figure~\ref{fig:uniform} plots the $200$ simulated signal pairs, while the top right panel gives the rotated observations.  The bottom left panel plots the recovered signal using the proposed nonparametric maximum likelihood method. Also included in the bottom right panel of the figure are the estimated marginal densities of the two sources of signal.

\begin{figure}[htbp]
\begin{center}
\includegraphics[height=0.7\textheight,width=0.5\textwidth,angle=270]{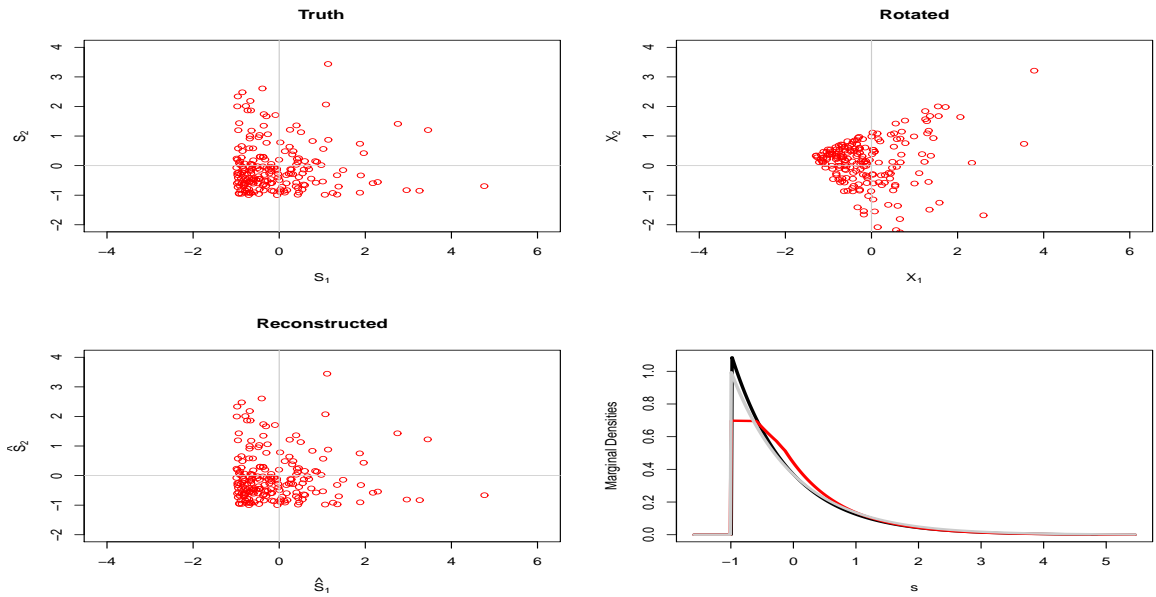}
\caption{Exponential signal: Top left panel, top right panel and bottom left panel give the true signal, rotated observations and the reconstructed signal respectively. The bottom right panel gives the estimated marginal densities along with the true marginal (grey line).}
\label{fig:exp}
\end{center}
\end{figure}

Figure~\ref{fig:exp} gives corresponding plots when the marginals have an $\mathrm{Exp}(1)-1$ distribution.  We note that both uniform and exponential distributions have log-concave densities and therefore our method not only recovers the mixing matrix but also accurately estimates the marginal densities, as can be seen in Figures~\ref{fig:uniform} and~\ref{fig:exp}.

\begin{figure}[htbp]
\begin{center}
\includegraphics[height=0.7\textheight,width=0.5\textwidth,angle=270]{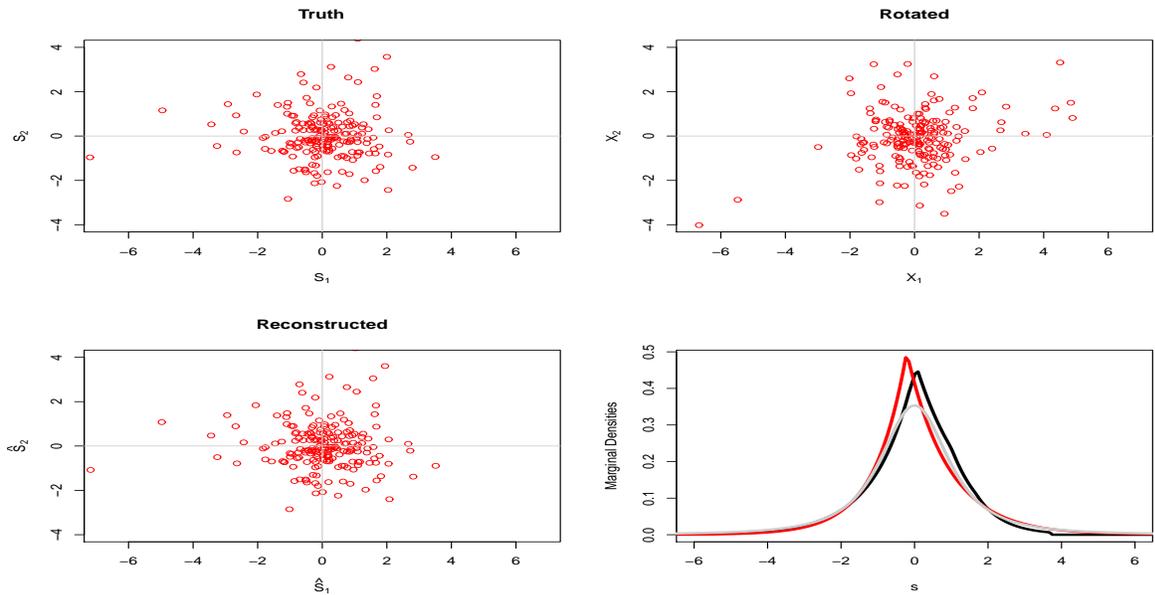}
\caption{$t_2$ signal: Top left panel, top right panel and bottom left panel give the true signal, rotated observations and the reconstructed signal respectively. The bottom right panel gives the estimated marginal densities along with the true marginal (grey line).}
\label{fig:t2}
\end{center}
\end{figure}

\begin{figure}[htbp]
\begin{center}
\includegraphics[height=0.7\textheight,width=0.5\textwidth,angle=270]{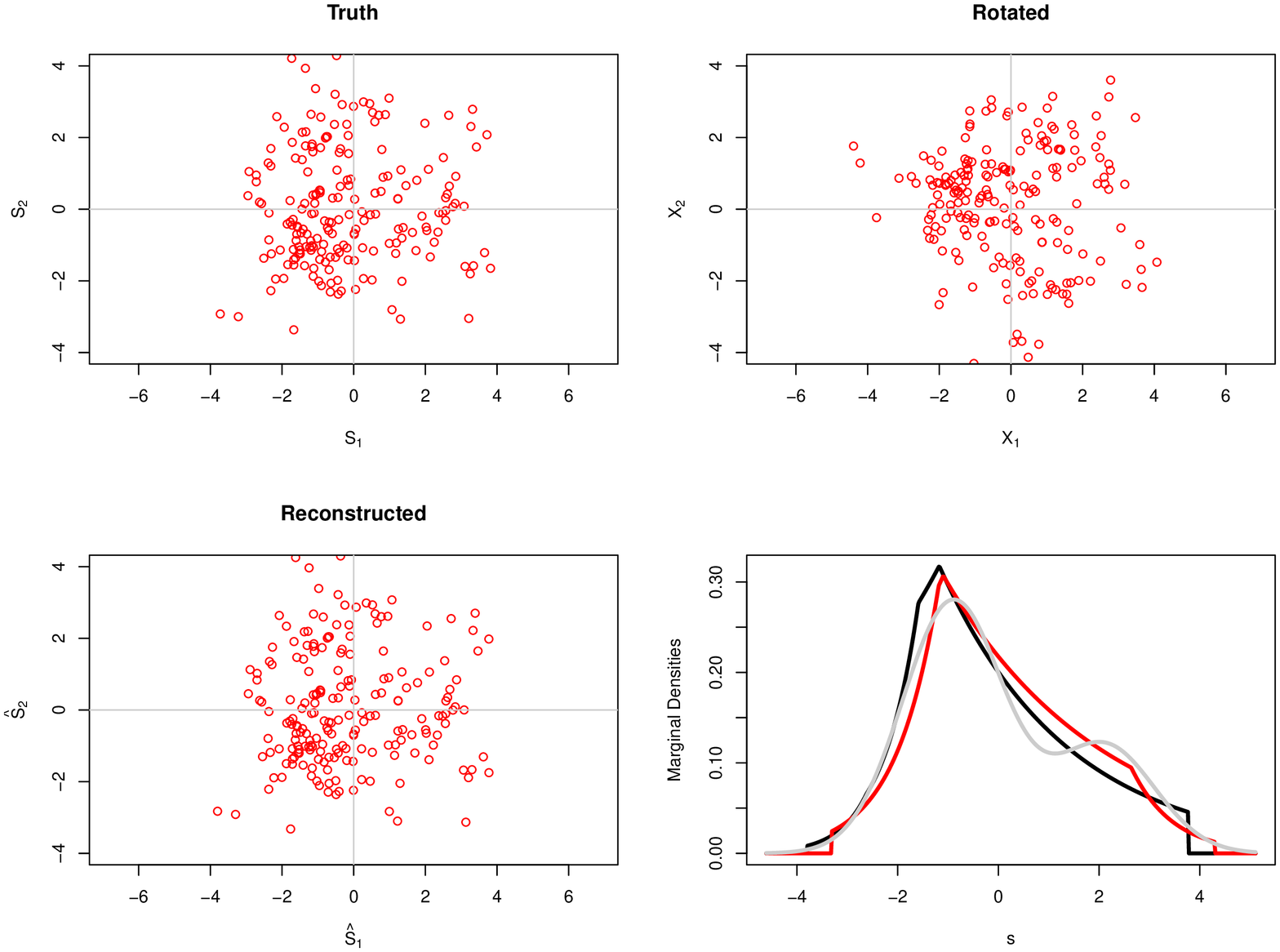}
\caption{Mixture of normals signal: Top left panel, top right panel and bottom left panel give the true signal, rotated observations and the reconstructed signal respectively. The bottom right panel gives the estimated marginal densities along with the true marginal (grey line).}
\label{fig:mix}
\end{center}
\end{figure}

To investigate the robustness of the proposed method when the marginal components do not have log-concave densities, we repeated the simulation in two other cases, with the true signal simulated firstly from a $t$-distribution with two degrees of freedom scaled by a factor of $1/\sqrt{2}$ and secondly from a mixture of normals distribution $0.7 N(-0.9,1)+0.3N(2.1,1)$.  Figures~\ref{fig:t2} and~\ref{fig:mix} show that, in both cases, the misspecification of the marginals does not affect the recovery of the signal.  Also, the estimated marginals represent estimates of the log-concave projection of the true marginals (a standard Laplace density in this case), as correctly predicted by our theoretical results.

\begin{figure}[htbp]
\begin{center}
\includegraphics[height=0.7\textheight,width=0.5\textwidth,angle=270]{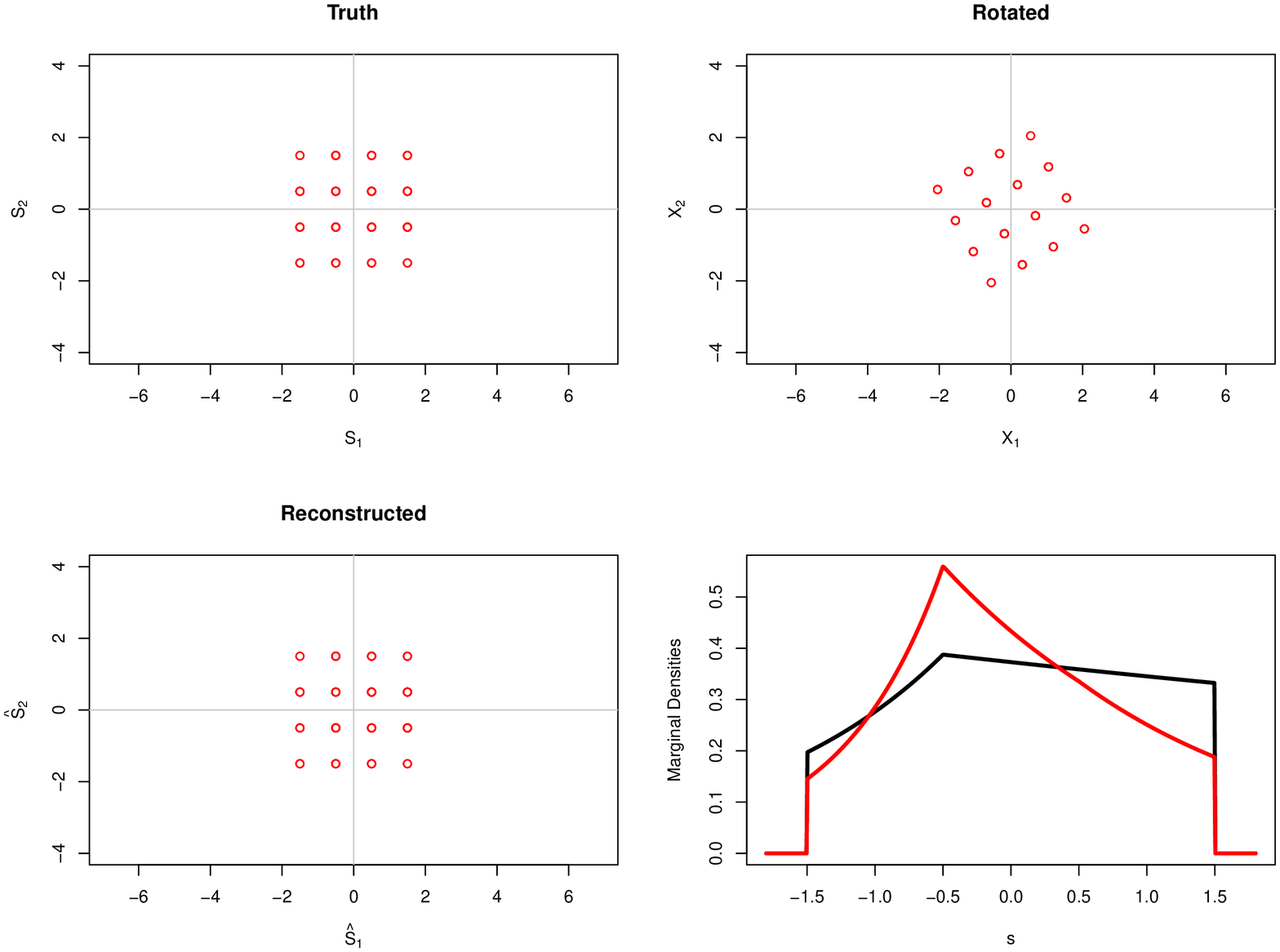}
\caption{Binomial signal: Top left panel, top right panel and bottom left panel give the true signal, rotated observations and the reconstructed signal respectively. The bottom right panel gives the estimated marginal densities.}
\label{fig:bin}
\end{center}
\end{figure}

As discussed before, one of the unique advantages of the proposed method over existing ones is its general applicability.  For example, the method can be used even when the marginal distributions of the true signal do not have densities.  To demonstrate this property, we now consider simulating signals from a $\mathrm{Bin}(3,1/2)-1.5$ distribution.  To the best of our knowledge, none of the existing ICA methods are applicable for these types of problems.  The simulation results presented in Figure~\ref{fig:bin} suggest that the method works very well in this case.

\begin{figure}[htbp]
\begin{center}
\includegraphics[height=0.7\textheight,width=0.5\textwidth,angle=270]{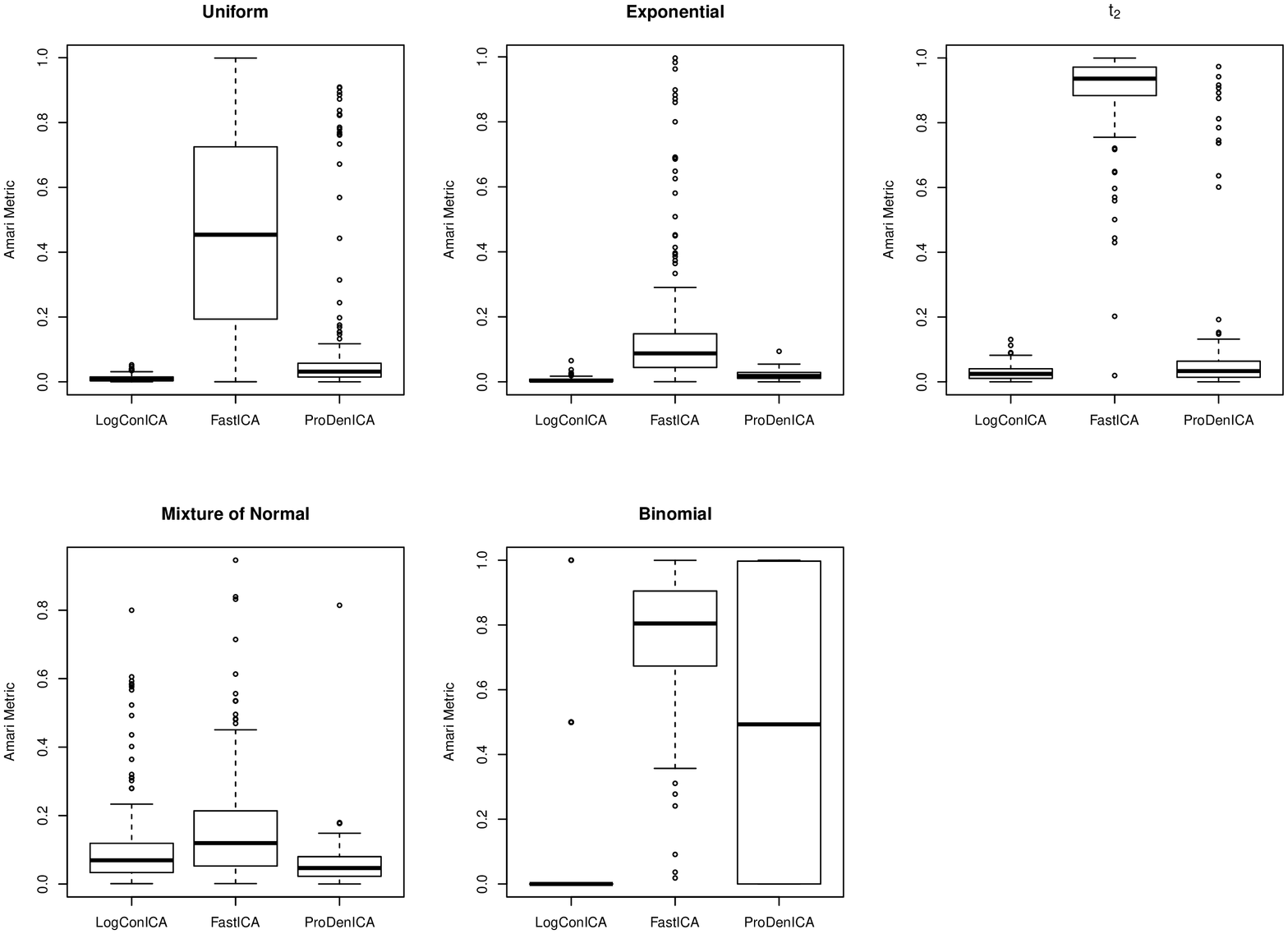}
\caption{Comparison between LogConICA, FastICA and ProDenICA.}
\label{fig:comp}
\end{center}
\end{figure}

To further conduct a comparative study, we repeated each of the previous simulations 200 times and computed our estimate along with those produced by the FastICA and ProDenICA methods.  FastICA is a popular parametric ICA method;  ProDenICA is a nonparametric ICA method proposed by \cite{HastieTibshirani2003}, and has been shown to enjoy the best performance among a large collection of existing ICA methods \citep{HTF2009}.  Both the FastICA and ProDenICA methods were implemented using the \verb+R+ package \verb+ProDenICA+ \citep{HastieTibshirani2010}.  To compare the performance of these methods, we follow convention \citep{HKO2001} and compute the Amari metric between the true unmixing matrix $W$ and its estimates.  The Amari metric between two $d\times d$ matrices is defined as
$$
M(A,B)=\frac{1}{2d}\sum_{i=1}^d\left(\frac{\sum_{j=1}^d|C_{ij}|}{\max_{1\le j\le d}|C_{ij}|}-1\right)+\frac{1}{2d}\sum_{j=1}^d\left(\frac{\sum_{i=1}^d|C_{ij}|}{\max_{1\le i\le d}|C_{ij}|}-1\right),
$$
where $C=(C_{ij})_{1\le i,j\le d}=AB^{-1}$. Boxplots of the Amari metric for all three methods are given in Figure~\ref{fig:comp}.

It is clear that both nonparametric methods outperform the parametric method.  Several further observations can also be made on the comparison between the two nonparametric methods.  For both uniform and exponential marginals, the proposed method improves upon ProDenICA. This might be expected since both distributions have log-concave densities. It is, however, interesting to note the robustness of the proposed method on the marginals as it still outperforms ProDenICA for $t_2$ marginals, and remains competitive for the mixture of normal marginals. The most significant advantage of the proposed method, however, is displayed when the marginals are binomial. Recall that ProDenICA, and perhaps all existing nonparametric methods, assume that the log density (or density itself) is smooth. This assumption is not satisfied with the binomial distribution and as a result, ProDenICA performs rather poorly.  In contrast, our proposed method works fairly well in this setting even though the true marginal does not have a log-concave density with respect to Lebesgue measure.  All these observations confirm our earlier theoretical development.


%

\section{Proofs}

\begin{prooftitle}{of Proposition~\ref{Prop:Cases}}
1. Suppose that $\int_{\mathbb{R}^d} \|x\| \, dP(x) = \infty$.  Fix an arbitrary $f \in \mathcal{F}_d^{\mathrm{ICA}}$, and find $\alpha > 0$ and $\beta \in \mathbb{R}$ such that $f(x) \leq e^{-\alpha\|x\| + \beta}$.  Then 
\[
\int_{\mathbb{R}^d} \log f \, dP \leq -\alpha \int_{\mathbb{R}^d} \|x\| \, dP(x) + \beta = -\infty.
\]
Thus $L^{**}(P) = -\infty$ and $\psi^{**}(P) = \mathcal{F}_d^{\mathrm{ICA}}$.

2. Now suppose that $\int_{\mathbb{R}^d} \|x\| \, dP(x) < \infty$, but $P(H) = 1$ for some hyperplane $H = \{x \in \mathbb{R}^d: a_1^\top x = \alpha\}$, where $a_1$ is a unit vector in $\mathbb{R}^d$ and $\alpha \in \mathbb{R}$.  Find $a_2,\ldots,a_d$ such that $a_1,\ldots,a_d$ is an orthonormal basis for $\mathbb{R}^d$.  Define the family of density functions
\[
f_\sigma(x) = \frac{1}{2\sigma} e^{-|a_1^{\sf T} x - \alpha|/\sigma}\prod_{j=2}^d \frac{e^{-|a_j^{\sf T} x|}}{2}.
\]
Then $f_\sigma \in \mathcal{F}_d^{\mathrm{ICA}}$, and
\begin{align*}
\int_{\mathbb{R}^d} \log f_\sigma(x) \, dP(x) &= -\log(\sigma) - d \log 2 - \sum_{j=2}^d \int_H |a_j^{\sf T} x| \, dP(x) \\
&\geq -\log(\sigma) - d \log 2 - \sum_{j=2}^d \int_H \|x\| \, dP(x) \rightarrow \infty
\end{align*}
as $\sigma \rightarrow 0$.

3. Now suppose that $P \in \mathcal{P}_d$.  Notice that the density $f(x) = 2^{-d} \prod_{j=1}^d e^{-|x_j|}$ belongs to $\mathcal{F}_d^{\mathrm{ICA}}$ and satisfies 
\[
\int_{\mathbb{R}^d} \log f \, dP = - \sum_{j=1}^d \int_{\mathbb{R}^d} |x_j| \, dP(x) - d \log 2 > -\infty.
\]
Moreover,
\[
\sup_{f \in \mathcal{F}_d^{\mathrm{ICA}}} \int_{\mathbb{R}^d} \log f \, dP \leq \sup_{f \in \mathcal{F}_d} \int_{\mathbb{R}^d} \log f \, dP < \infty,
\]
where the second inequality follows from the proof of Theorem~2.2 of \citet{DSS2011}.  We may therefore take a sequence $f^1,f^2,\ldots \in \mathcal{F}_d^{\mathrm{ICA}}$ such that 
\[
\int_{\mathbb{R}^d} \log f^n \, dP \nearrow \sup_{f \in \mathcal{F}_d^{\mathrm{ICA}}} \int_{\mathbb{R}^d} \log f \, dP.
\]
Let $\mathrm{csupp}(P)$ denote the convex support of $P$; that is, the intersection of all closed, convex sets having $P$-measure 1.  Following the arguments in the proof of Theorem~2.2 of \citet{DSS2011}, there exist $\alpha > 0$ and $\beta \in \mathbb{R}$ such that $\sup_{n \in \mathbb{N}} f^n(x) \leq e^{-\alpha \|x\| + \beta}$ for all $x \in \mathbb{R}^d$.  Moreover, these arguments (see also the proof of Theorem~4 of \citet{CuleSamworth2010}) yield the existence of a closed, convex set $C \supseteq \mathrm{int}(\mathrm{csupp}(P))$, a log-concave density $f^{**} \in \mathcal{F}_d$ with $\{x \in \mathbb{R}^d: f^{**}(x) > 0\} = C$ and a subsequence $(f^{n_k})$ such that 
\[
f^{**}(x) = \lim_{k \rightarrow \infty} f^{n_k}(x) \quad \text{for all $x \in \mathrm{int}(C) \cup (\mathbb{R}^d \setminus C)$}.
\]
Since the boundary of $C$ has zero Lebesgue measure, we deduce from Fatou's lemma applied to the non-negative functions $x \mapsto e^{-\alpha \|x\| + \beta} - f^{n_k}(x)$ that
\[
\int_{\mathbb{R}^d} \log f^{**} \, dP \geq \limsup_{k \rightarrow \infty} \int_{\mathbb{R}^d} \log f^{n_k} \, dP = \sup_{f \in \mathcal{F}_d^{\mathrm{ICA}}} \int_{\mathbb{R}^d} \log f \, dP.
\]
It remains to show that $f^{**} \in \mathcal{F}_d^{\mathrm{ICA}}$.  We can write
\[
f^{n_k}(x) = |\det W^k| \prod_{j=1}^d f_j^k((w_j^k)^{\sf T} x),
\]
where $W^k \in \mathcal{W}$ and $f_j^k \in \mathcal{F}_1$ for each $k \in \mathbb{N}$ and $j=1,\ldots,d$.  Let $X^k$ be a random vector with density $f^{n_k} \in \mathcal{F}_d^{\mathrm{ICA}}$, and let $X$ be a random vector with density $f^{**} \in \mathcal{F}_d$.  We know that $X^k \stackrel{d}{\rightarrow} X$ as $k \rightarrow \infty$, and that $(w_1^k)^{\sf T}X^k,\ldots,(w_d^k)^{\sf T}X^k$ are independent for each $k$.  Let $\tilde{w}_j^k = w_j^k/\|w_j^k\|$ and $\tilde{f}_j^k(x) = \|w_j^k\|f_j^k(\|w_j^k\|x)$.  Then we have
\begin{equation}
\label{Eq:SecondRep}
f^{n_k}(x) = |\det \tilde{W}^k| \prod_{j=1}^d \tilde{f}_j^k((\tilde{w}_j^k)^{\sf T} x),
\end{equation}
where the matrix $\tilde{W}^k$ has $j$th row $\tilde{w}_j^k$.  Moreover, $\tilde{W}^k \in \mathcal{W}$ and $\tilde{f}_1^k,\ldots,\tilde{f}_d^k \in \Fcal_1$, so~(\ref{Eq:SecondRep}) provides an alternative, equivalent representation of the density $f^{n_k}$, in which each row of the unmixing matrix has unit Euclidean length.  By reducing to a further subsequence if necessary, we may assume that for each $j=1,\ldots,d$, there exists $\tilde{w}_j \in \mathbb{R}^d$ such that $\tilde{w}_j^k \rightarrow \tilde{w}_j$ as $k \rightarrow \infty$.  By Slutsky's theorem, it then follows that 
\[
((\tilde{w}_1^k)^{\sf T} X^k,\ldots,(\tilde{w}_d^k)^{\sf T} X^k) \stackrel{d}{\rightarrow} (\tilde{w}_1^{\sf T} X,\ldots,\tilde{w}_d^{\sf T} X).
\]
Thus, for any $t \in \mathbb{R}^d$,
\begin{align*}
\mathbb{E}(e^{it^{\sf T}(\tilde{w}_1^{\sf T} X,\ldots,\tilde{w}_d^{\sf T} X)}) &= \lim_{k \rightarrow \infty} \mathbb{E}(e^{it^{\sf T}((\tilde{w}_1^k)^{\sf T} X^k,\ldots,(\tilde{w}_d^k)^{\sf T} X^k)}) \\
&= \lim_{k \rightarrow \infty} \prod_{j=1}^d \mathbb{E}(e^{it_j (\tilde{w}_j^k)^{\sf T} X^k}) = \prod_{j=1}^d \mathbb{E}(e^{it_j \tilde{w}_j^{\sf T} X}).
\end{align*}
We conclude that $\tilde{w}_1^{\sf T} X,\ldots,\tilde{w}_d^{\sf T} X$ are independent.  Since $\|\tilde{w}_j\| = 1$ for all $j$, we deduce further that $\tilde{W} = (\tilde{w}_1,\ldots,\tilde{w}_d)^{\sf T}$ is non-singular.  Moreover, each of $\tilde{w}_1^{\sf T} X,\ldots,\tilde{w}_d^{\sf T} X$ has a log-concave density, by Theorem~6 of \citet{Prekopa1973}.  This shows that $f^{**} \in \mathcal{F}_d^{\mathrm{ICA}}$, as required.
\hfill $\Box$
\end{prooftitle}

\begin{prooftitle}{of Theorem~\ref{Thm:PdICA}}
Suppose that $P \in \mathcal{P}_d^{\mathrm{ICA}}$ satisfies
\[
P(B) = \prod_{j=1}^d P_j(w_j^{\sf T} B)
\]
for some $W \in \mathcal{W}$ and $P_1,\ldots,P_d \in \mathcal{P}_1$.  Consider maximising
\[
\int_{\mathbb{R}^d} \log f(x) \, dP(x)
\]
over $f \in \mathcal{F}_d$.  Letting $s = Wx$ and $\tilde{f}(s) = f(As)$, where $A = W^{-1}$, we can equivalently maximise
\[
\int_{\mathbb{R}^d} \log \tilde{f}(s) \, d\biggl(\bigotimes_{j=1}^d P_j(s_j)\biggr)
\]
over $\tilde{f} \in \mathcal{F}_d$.  But, by Theorem~4 of \citet{ChenSamworth2012}, the unique solution to this maximisation problem is to choose $\tilde{f}(z) = \prod_{j=1}^d f_j^*(z_j)$, where $f_j^* = \psi^*(P_j)$. This shows that $f^* := \psi^*(P)$ can be written as
\[
f^*(x) = |\det W| \prod_{j=1}^d f_j^*(w_j^{\sf T} x),
\]
Since $f^* \in \mathcal{F}_d^{\mathrm{ICA}}$ also, we deduce that $f^*$ is also the unique maximiser of $\int_{\mathbb{R}^d} \log f \, dP$ over $f \in  \mathcal{F}_d^{\mathrm{ICA}}$, so $\psi^{**}(P) = \psi^*(P)$.
\hfill $\Box$
\end{prooftitle}

\begin{prooftitle}{of Theorem~\ref{Thm:Identifiability}}
Suppose that $P \in \mathcal{P}_d^{\mathrm{ICA}}$.  Let $X \sim P$, so there exists $W \in \mathcal{W}$ such that $WX$ has independent components.  Writing $P_j$ for the marginal distribution of $w_j^{\sf T}X$, note that $P_1,\ldots,P_d \in \mathcal{P}_1$.  By Theorem~\ref{Thm:PdICA} and the identifiability result of \citet{ErikssonKoivunen2004}, it therefore suffices to show that $P_j \in \mathcal{P}_1$ has a Gaussian density if and only if $\psi^*(P_j)$ is a Gaussian density.  If $P_j$ has a Gaussian density $f_j^*$, then since $f_j^*$ is log-concave, we have $f_j^* = \psi^*(P_j)$.  Conversely, suppose that $P_j$ does not have a Gaussian density.  Since $f_j^* = \psi^*(P_j)$ satisfies $\int_{-\infty}^\infty x \, dP_j(x) = \int_{-\infty}^\infty x f_j^*(x) \, dx$ \citep[][Remark~2.3]{DSS2011}, we may assume without loss of generality that $P_j$ and $f_j^*$ have mean zero.  We consider maximising 
\[
\int_{-\infty}^\infty \log f \, dP_j
\]
over all mean zero Gaussian densities $f$.  Writing $\phi_{\sigma^2}$ for the mean zero Gaussian density with variance $\sigma^2$, we have
\[
\int_{-\infty}^\infty \log \phi_{\sigma^2} \, dP_j = -\frac{1}{2\sigma^2} \int_{-\infty}^\infty x^2 \, dP_j(x) - \frac{1}{2} \log (2\pi\sigma^2).
\]
This expression is maximised uniquely in $\sigma^2$ at $\sigma_*^2 = \int_{-\infty}^\infty x^2 \, dP_j(x)$.  But \citet{ChenSamworth2012} show that the only way a distribution $P_j$ and its log-concave projection $\psi^*(P_j)$ can have the same second moment is if $P_j$ has a log-concave density, in which case $P_j$ has density $\psi^*(P_j)$.  We therefore conclude that the only way $\psi^*(P_j)$ can be a Gaussian density is if $P_j$ has a Gaussian density, a contradiction.
\hfill $\Box$
\end{prooftitle}
\begin{prooftitle}{of Proposition~\ref{Prop:Pd}}
The proof of this proposition is very similar to the proof of Theorem~4.5 of \citet{DSS2011}, so we only sketch the argument here.  For each $n \in \mathbb{N}$, let $f^n \in \psi^{**}(P^n)$, and consider an arbitrary subsequence $(f^{n_k})$.  By reducing to a further subsequence if necessary, we may assume that $L^{**}(P^{n_k}) \rightarrow \lambda \in [-\infty,\infty]$.  Observe that 
\[
\lambda \geq \lim_{k \rightarrow \infty} \int_{\mathbb{R}^d} \log\bigl(2^{-d}e^{-\sum_{j=1}^d |x_j|}\bigr) \, dP^{n_k}(x) = -d \log 2 - \sum_{j=1}^d \int_{\mathbb{R}^d} |x_j| \, dP(x) > -\infty.
\]
Arguments from convex analysis can be used to show that the sequence $(f^{n_k})$ is uniformly bounded above, and $\liminf_{k \in \mathbb{N}} f^{n_k}(x_0) > -\infty$ for all $x_0 \in \mathrm{int}(\mathrm{csupp}(P))$.  From this it follows that there exist $a > 0$ and $b \in \mathbb{R}$ such that $\sup_{k \in \mathbb{N}} \sup_{x \in \mathbb{R}^d} f^{n_k}(x) \leq -a\|x\| + b$.  Thus, by reducing to a further subsequence if necessary, we may assume there exists $f^{**} \in \mathcal{F}_d$ such that 
\begin{align}
\label{Eq:aeconv}
\limsup_{k \rightarrow \infty, x \rightarrow x_0} f^{n_k}(x) &= f^{**}(x_0) \quad \text{for all $x_0 \in \mathbb{R}^d \setminus \partial \{x \in \mathbb{R}^d: f^{**}(x) > 0\}$} \\
 \limsup_{k \rightarrow \infty, x \rightarrow x_0} f^{n_k}(x) &\leq f^{**}(x_0) \quad \text{for all $x_0 \in \partial \{x \in \mathbb{R}^d: f^{**}(x) > 0\}.$} \nonumber
 \end{align}
 Note from this that
 \[
 \lambda = \lim_{k \rightarrow \infty} \int_{\mathbb{R}^d} \log f^{n_k} \, dP^{n_k} \leq -a\int_{\mathbb{R}^d} \|x\| \, dP(x) + b < \infty.
 \]
 In fact, we can use the argument from the proof of Proposition~\ref{Prop:Cases} to deduce that $f^{**} \in \mathcal{F}_d^{\mathrm{ICA}}$.  Skorokhod's representation theorem and Fatou's lemma can then be used to show that $\lambda \leq \int_{\mathbb{R}^d} \log f^{**} \, dP \leq L^{**}(P)$.
 
We can obtain the other bound $\lambda \geq L^{**}(P)$ by taking any element of $\psi^{**}(P)$, approximating it from above using Lipschitz continuous functions, as in the proof of Theorem~4.5 of \citet{DSS2011}, and using monotone convergence.  From these arguments, we conclude that $L^{**}(P^n) \rightarrow L^{**}(P)$ and $f^{**} \in \psi^{**}(P)$.

We can see from~(\ref{Eq:aeconv}) that $f^{n_k} \stackrel{a.e.}{\rightarrow} f^{**}$, so $\int_{\mathbb{R}^d} |f^{n_k} - f^{**}| \rightarrow 0$, by Scheff\'e's theorem.  Thus, given any $f^n \in \psi^{**}(P^n)$ and any subsequence $(f^{n_k})$, we can find $f^{**} \in \psi^{**}(P)$ and a further subsequence of $(f^{n_k})$ which converges to $f^{**}$ in total variation distance.  This yields the second part of the proposition. 
\hfill $\Box$
\end{prooftitle}

\begin{prooftitle}{of Theorem~\ref{Thm:Conv}}
The first part of the theorem is a special case of Proposition~\ref{Prop:Pd}.  Now suppose $P \in \mathcal{P}_d^{\mathrm{ICA}}$ is identifiable and is represented by $W \in \mathcal{W}$ and $P_1,\ldots,P_d \in \mathcal{P}_1$.  Suppose without loss of generality that $\|w_j\| = 1$ for all $j=1,\ldots,d$ and let $f^{**} = \psi^{**}(P)$.  Recall from Theorem~\ref{Thm:PdICA} that if $X$ has density $f^{**}$, then $w_j^{\sf T}X$ has density $f_j^* = \psi^*(P_j)$.

Suppose for a contradiction that we can find $\epsilon > 0$, integers $1 \leq n_1 < n_2 < \ldots$, $f^k \in \psi^{**}(P^{n_k})$ and $(W^k,f_1^k,\ldots,f_d^k) \stackrel{\mathrm{ICA}}{\sim} f^k$ such that 
\[
\inf_{k \in \mathbb{N}} \inf_{\epsilon_j^k \in \mathbb{R} \setminus \{0\}} \inf _{\pi^k \in \Pi_d} \biggl\{\|(\epsilon_j^k)^{-1} w_{\pi^k(j)}^k  - w_j\| + \int_{-\infty}^\infty \bigl| |\epsilon_j^k| f_{\pi^k(j)}^k(\epsilon_j^k x) - f_j^*(x)\bigr| \, dx\biggr\} \geq \epsilon. 
\]
We can find a subsequence $1 \leq k_1 < k_2 < \ldots$ such that $w_j^{k_l}/\|w_j^{k_l}\| \rightarrow \tilde{w}_j$, say, as $l \rightarrow \infty$, for all $j=1,\ldots,d$.    The argument towards the end of the proof of Case 3 of Proposition~\ref{Prop:Cases} shows that $\tilde{W}$ can be used to represent the unmixing matrix of $f^{**}$, so by the identifiability result of \citet{ErikssonKoivunen2004} and the fact that $\|\tilde{w}_j\| = 1$, there exist $\tilde{\epsilon}_1,\ldots,\tilde{\epsilon}_d \in \{-1,1\}$ and a permutation $\pi$ of $\{1,\ldots,d\}$ such that $\tilde{\epsilon}_j \tilde{w}_{\pi(j)} = w_j$.  Setting $\pi^n = \pi$ and $\epsilon_j^n = \tilde{\epsilon}_j^{-1}\|w_{\pi^n(j)}^n\|$, we deduce that
\[
(\epsilon_j^{k_l})^{-1} w_{\pi^{k_l}(j)}^{k_l} = \tilde{\epsilon}_j \frac{w_{\pi(j)}^{k_l}}{\|w_{\pi(j)}^{k_l}\|} \rightarrow w_j,
\]
for $j=1,\ldots,d$.  Now observe that if $X^{k_l}$ has density $f^{k_l}$, then by Slutsky's theorem, $(\epsilon_j^{k_l})^{-1} (w_{\pi^{k_l}(j)}^{k_l})^{\sf T}X^{k_l} \stackrel{d}{\rightarrow} w_j^{\sf T}X$.  It therefore follows from Proposition~2(c) of \citet{CuleSamworth2010} that
\[
\int_{-\infty}^\infty \bigl| |\epsilon_j^{k_l}| f_{\pi^{k_l}(j)}^n(\epsilon_j^{k_l} x) - f_j^*(x)\bigr| \, dx \rightarrow 0
\]
for $j=1,\ldots,d$.  This contradiction establishes that  
\begin{align}
\label{Eq:supinf}
\sup_{f^n \in \psi^{**}(P^n)} \sup_{(W^n,f_1^n,\ldots,f_d^n) \stackrel{\mathrm{ICA}}{\sim} f^n} \inf_{\pi^n \in \Pi_d} \inf_{\epsilon_1^n,\ldots,\epsilon_d^n \in \mathbb{R} \setminus \{0\}} \biggl\{\|&(\epsilon_j^n)^{-1} w_{\pi^n(j)}^n - w_j\| \nonumber \\
&+ \int_{-\infty}^\infty \bigl| |\epsilon_j^n| f_{\pi^n(j)}^n(\epsilon_j^n x) - f_j^*(x)\bigr| \, dx\biggr\} \rightarrow 0,
\end{align}
for each $j=1,\ldots,d$.

It remains to prove that for sufficiently large $n$, every $f^n \in \psi^{**}(P^n)$ is identifiable.  Recall from the identifiability result of \citet{ErikssonKoivunen2004} and Theorem~\ref{Thm:Identifiability} that not more than one of $f_1^*,\ldots,f_d^*$ is Gaussian.  Let $\phi_{\mu,\sigma^2}(\cdot)$ denote the univariate normal density with mean $\mu$ and variance $\sigma^2$.  Let $J$ denote the index set of the non-Gaussian densities among $f_1^*,\ldots,f_d^*$, so the cardinality of $J$ is at least $d-1$, and consider, for each $j \in J$, the problem of minimising $g(\mu,\sigma) = \int_{-\infty}^\infty |\phi_{\mu,\sigma^2} - f_j^*|$ over $\mu \in \mathbb{R}$ and $\sigma > 0$.  Observe that $g$ is continuous with $g(\mu,\sigma) < 2$ for all $\mu$ and $\sigma$, that $\inf_{\mu \in \mathbb{R}} g(\mu,\sigma) \rightarrow 2$ as $\sigma \rightarrow 0,\infty$ and $\inf_{\sigma > 0} g(\mu,\sigma) \rightarrow 2$ as $|\mu| \rightarrow \infty$.  It follows that $g$ attains its infimum, and there exists $\eta > 0$ such that
\begin{equation}
\label{Eq:Inf}
\inf_{j \in J} \inf_{\mu \in \mathbb{R}} \inf_{\sigma > 0} \int_{-\infty}^\infty |\phi_{\mu,\sigma^2} - f_j^*| \geq \eta.
\end{equation}
Comparing~(\ref{Eq:supinf}) and~(\ref{Eq:Inf}), we see that, for sufficiently large $n$, whenever $f^n \in \psi^{**}(P^n)$ and $(W^n,f_1^n,\ldots,f_d^n) \stackrel{\mathrm{ICA}}{\sim} f^n$, at most one of the densities $f_1^n,\ldots,f_d^n$ can be Gaussian.  It follows that when $n$ is large, every $f^n \in \psi^{**}(P^n)$ is identifiable.
\hfill $\Box$
\end{prooftitle}

\begin{prooftitle}{of Proposition~\ref{Prop:EmpLike}}
It is well-known that for fixed $W \in \mathcal{W}$, the nonparametric likelihood $L(\cdot)$ defined in~(\ref{Eq:NonLike}) is maximised by choosing
\[
\hat{P}_j^W = \frac{1}{n}\sum_{i=1}^n \delta_{w_j^{\sf T}\bx_i}, \quad j=1,\ldots,d.
\]
For $i = 1\ldots,n$, $W \in \mathcal{W}$ and $j=1,\ldots,d$, let 
\[
n_{w_j}(i) = \bigl\{\tilde{i} \in \{1,\ldots,n\}: w_j^{\sf T}\bx_{\tilde{i}} = w_j^{\sf T}\bx_i\bigr\}.
\]
The binary relation $i \sim \tilde{i}$ if $n_{w_j}(i) = n_{w_j}(\tilde{i})$ defines an equivalence relation on $\{1,\ldots,n\}$, so we can let $I^W$ denote a set of indices obtained by choosing one element from each equivalence class.  Then
\[
L(W,\hat{P}_1^W,\ldots,\hat{P}_d^W) = \prod_{j=1}^d \frac{|n_{w_j}(1)||n_{w_j}(2)|\ldots |n_{w_j}(n)|}{n^n} = \prod_{j=1}^d n^{-n} \prod_{i \in I^W} |n_{w_j}(i)|^{|n_{w_j}(i)|}
\]
Since $\bx_1,\ldots,\bx_n$ are in general position by hypothesis, we have that $\sum_{i \in I^W} (|n_{w_j}(i)| - 1) \leq d - 1$.  It follows that $L(W,\hat{P}_1^W,\ldots,\hat{P}_d^W) \leq (d^d/n^n)^d$.  Moreover, for any choice $J$ of distinct indices in $\{1,\ldots,n\}$ if we construct the matrix $W_J \in \mathcal{W}$ as described just before the statement of Proposition~\ref{Prop:EmpLike}, then $L(W_J,\hat{P}_1^{W_J},\ldots,\hat{P}_d^{W_J}) = (d^d/n^n)^d$. 
\hfill $\Box$
\end{prooftitle}

\begin{prooftitle}{of Corollary~\ref{Cor:Hats}}
Let $\hat{P}^{n,\bz}$ denote the empirical distribution of $\bz_1 = \hat{\Sigma}^{-1/2}\bx_1,\ldots,\bz_n = \hat{\Sigma}^{-1/2}\bx_n$.  Writing $\bar{\bz} = n^{-1}\sum_{i=1}^n \bz_i$ and $\bar{\bx} = n^{-1}\sum_{i=1}^n \bx_i$, note that the covariance matrix corresponding to $\hat{P}^{n,\bz}$ is
\[
\frac{1}{n}\sum_{i=1}^n (\bz_i - \bar{\bz})(\bz_i - \bz)^{\sf T} = \frac{1}{n}\sum_{i=1}^n \hat{\Sigma}^{-1/2} (\bx_i - \bar{\bx})(\bx_i - \bx)^{\sf T}\hat{\Sigma}^{-1/2} = I.
\]
Observe further that there is a bijection between the set of maximisers $(\hat{\hat{W}}^n,\hat{\hat{f}}_1^n,\ldots,\hat{\hat{f}}_d^n)$ of $\ell^n(W,f_1,\ldots,f_d;\bx_1,\ldots,\bx_n)$ over $W \in O(d)\hat{\Sigma}^{-1/2}$ and $f_1,\ldots,f_d \in \Fcal_1$, and the set of maximisers $(\hat{O}^n,\hat{g}_1^n,\ldots,\hat{g}_d^n)$ of $\ell^n(O,g_1,\ldots,g_d;\bz_1,\ldots,\bz_n)$ over $O \in O(d)$ and $g_1,\ldots,g_d \in \mathcal{F}_1$ via the correspondence $\hat{\hat{W}}^n = \hat{O}^n\hat{\Sigma}^{-1/2}$ and $\hat{\hat{f}}_j^n = \hat{g}_j^n$.

It follows from the discussion in Section~\ref{Sec:Pre-Whiten} that maximising $\ell^n(O,g_1,\ldots,g_d;\bz_1,\ldots,\bz_n)$ over $O \in O(d)$ and $g_1,\ldots,g_d \in \mathcal{F}_1$ amounts to computing the log-concave ICA projection of $\hat{P}^{n,\bz}$.  Existence of a maximiser therefore follows from Proposition~\ref{Prop:Cases} and the fact that the convex hull of $\bz_1,\ldots,\bz_n$ is $d$-dimensional with probability 1 for sufficiently large~$n$.

Now suppose $\hat{O}^n$ and $\hat{g}_1^n,\ldots,\hat{g}_d^n$ represent the log-concave ICA projection $\psi^{**}(\hat{P}^{n,\bz})$.  Further, let $P^{0,\bz}$ denote the distribution of $\Sigma^{-1/2}\bx_1$, so $P^{0,\bz} \in \mathcal{P}_d^{\mathrm{ICA}}$ has identity covariance matrix and suppose $(O^0,P_1^{0,\bz},\ldots,P_d^{0,\bz}) \stackrel{\mathrm{ICA}}{\sim} P^{0,\bz}$.  Then $d(\hat{P}^{n,\bz},P^{0,\bz}) \stackrel{a.s.}{\rightarrow} 0$ as $n \rightarrow \infty$, so by Theorem~\ref{Thm:Conv}, there exist a permutation $\hat{\hat{\pi}}^n$ of $\{1,\ldots,d\}$ and scaling factors $\hat{\hat{\epsilon}}_1^n,\ldots,\hat{\hat{\epsilon}}_d^n \in \mathbb{R}\setminus  \{0\}$ such that
\[
(\hat{\hat{\epsilon}}_j^n)^{-1}\hat{o}_{\hat{\hat{\pi}}^n(j)}^n \stackrel{a.s.}{\rightarrow} o_j^0 \quad \text{and} \quad \int_{-\infty}^\infty \bigl| |\hat{\hat{\epsilon}}_j^n|\hat{g}_{\hat{\hat{\pi}}^n(j)}^n(\hat{\hat{\epsilon}}_j^n x) - g_j^*(x)\bigr| \, dx \stackrel{a.s.}{\rightarrow} 0,
\]
where $g_j^* = \psi^*(P_j^{0,\bz})$.  Writing $W^0 = O^0\Sigma^{-1/2}$, $\hat{\hat{W}}^n = \hat{O}^n\hat{\Sigma}^{-1/2}$, $\hat{\hat{f}}_j^n = \hat{g}_j^n$ and noting that $g_j^* = \psi^*(P_j^{0,\bz}) = \psi^*(P_j^0) = f_j^*$, the conclusion of the corollary follows immediately.
\hfill $\Box$
\end{prooftitle}

\begin{prooftitle}{of Proposition~\ref{Prop:Diff}}
For $\epsilon > 0$, let $W_\epsilon = W\exp(\epsilon Y)$, and let $w_{j,\epsilon}$ denote the $j$th row of $W_\epsilon$.  Notice that
\[
w_{j,\epsilon}^{\sf T}\bx_i = w_j^{\sf T}\bx_i + \epsilon c_j^{\sf T}\bx_i + O(\epsilon^2)
\]
as $\epsilon \searrow 0$.  It follows that for sufficiently small $\epsilon > 0$,
\begin{align*}
\frac{g(W_\epsilon) - g(W)}{\epsilon} &= \frac{1}{\epsilon}\sum_{i=1}^n \sum_{j=1}^d \Bigl\{\min_{k=1,\ldots,m_j} (b_{jk} w_{j,\epsilon}^{\sf T}\bx_i - \beta_{jk}) - \min_{k=1,\ldots,m_j} (b_{jk} w_j^{\sf T}\bx_i - \beta_{jk})\Bigr\} \\
&= \frac{1}{\epsilon}\sum_{i=1}^n \sum_{j=1}^d b_{jk_{ij}} (w_{j,\epsilon}^{\sf T} \bx_i - w_j^{\sf T}\bx_i) \\
&\rightarrow \sum_{i=1}^n \sum_{j=1}^d b_{jk_{ij}} c_j^{\sf T} \bx_i
\end{align*}
as $\epsilon \searrow 0$.
\hfill $\Box$ 
\end{prooftitle}


\begin{thebibliography}{99}

\bibitem[{Bach and Jordan(2002)}]{BachJordan2002}
Bach, F., Jordan, M. I. (2002) Kernel independent component analysis.
\newblock \emph{Journal of Machine Learning Research}, \textbf{3}, 1-48.

\bibitem[{Chen and Bickel(2005)}]{ChenBickel2005}
Chen, A. and Bickel, P. J. (2005) Consistent independent component analysis and pre-whitening.
\newblock \emph{IEEE Trans. Signal. Proc.}, \textbf{53}, 3625--3632.

\bibitem[{Chen and Bickel(2006)}]{ChenBickel2006}
Chen, A. and Bickel, P. J. (2006) Efficient independent component analysis, 
\newblock \emph{The Annals of Statistics}, \textbf{34}, 2825-2855.

\bibitem[{Chen and Samworth(2012)}]{ChenSamworth2012}
Chen, Y. and Samworth, R. J. (2012) Smoothed log-concave maximum likelihood estimation with applications.
\newblock Preprint, available at \verb+http://arxiv.org/pdf/1102.1191v4+.

\bibitem[{Comon(1994)}]{Comon1994}
Comon, P. (1994) Independent component analysis, A new concept?
\newblock \emph{Signal Proc.}, \textbf{36}, 287--314.

\bibitem[{Cule and Samworth(2010)}]{CuleSamworth2010}
Cule, M. and Samworth, R. (2010) Theoretical properties of the log-concave maximum likelihood estimator of a multidimensional density.
\newblock \emph{Elect. J. Statist.}, \textbf{4}, 254--270.

\bibitem[{Cule, Samworth and Stewart(2010)}]{CSS2010}
Cule, M, Samworth, R. and Stewart, M. (2010) Maximum likelihood estimation of a multi-dimensional log-concave density
\newblock \emph{J. Roy. Statist. Soc., Ser. B (with discussion)}, \textbf{72}, 545--607.

\bibitem[{D\"umbgen and Rufibach(2011)}]{DuembgenRufibach2011}
D\"umbgen, L. and Rufibach, K. (2011) \verb+logcondens+: Computations Related to Univariate Log-Concave Density Estimation. 
\newblock \emph{J. Statist. Software}, \textbf{39}, 1--28.

\bibitem[{D\"umbgen, Samworth and Schuhmacher(2011)}]{DSS2011}
D\"umbgen, L., Samworth, R. and Schuhmacher, D. (2011) Approximation by log-concave distributions, with applications to regression.
\newblock \emph{Ann. Statist.}, \textbf{39}, 702--730.

\bibitem[{Eriksson and Koivunen(2004)}]{ErikssonKoivunen2004}
Eriksson, J. and Koivunen, V. (2004) Identifiability, separability and uniqueness of linear ICA models.
\newblock \emph{IEEE Signal Processing Letters}, \textbf{11}, 601--604.

\bibitem[{Hastie and Tibshirani(2003)}]{HastieTibshirani2003}
Hastie, T. and Tibshirani, R. (2003) Independent component analysis through product density estimation. 
in \newblock \emph{Advances in Neural Information Processing Systems 15 (Becker, S. and Obermayer, K., eds)}, MIT Press, Cambridge, MA. pp 649-656.

\bibitem[{Hastie and Tibshirani(2010)}]{HastieTibshirani2010}
Hastie, T. and Tibshirani, R. (2003) \verb+ProDenICA+: Product Density Estimation for ICA using tilted Gaussian density estimates
\newblock \texttt{R} package version 1.0
\newblock \texttt{http://cran.r-project.org/web/packages/ProDenICA/}.


\bibitem[{Hastie, Tibshirani and Friedman(2009)}]{HTF2009}
Hastie, T., Tibshirani, R. and Friedman (2009) \newblock \emph{The Elements of Statistical Learning}, New York: Springer.

\bibitem[{Hyvarinen, Karhunen and Oja(2001)}]{HKO2001}
Hyv\"arinen, A., Karhunen, J. and Oja, E. (2001)
\newblock \emph{Independent Component Analysis}, New York: John Wiley \& Sons.

\bibitem[{Owen(1990)}]{Owen1990}Owen, A. (1990)
\newblock \emph{Empirical Likelihood}, Chapman and Hall, London.

\bibitem[{Pr\'ekopa(1973)}]{Prekopa1973}
Pr\'ekopa, A. (1973) On logarithmically concave measures and functions.
\newblock \emph{Acta Scientarium Mathematicarum}, \textbf{34}, 335--343.

\bibitem[{Rufibach and D\"umbgen(2006)}]{RufibachDuembgen2006} Rufibach, K. and D\"umbgen, L. (2006) \verb+logcondens+: Estimate a log-concave probability density from \emph{i.i.d} Observations
\newblock \texttt{R} package version 2.01
\newblock \texttt{http://cran.r-project.org/web/packages/logcondens/}.

\bibitem[{Samarov and Tsybakov(2004)}]{SamTsy2004} Samarov, A. and Tsybakov, A. (2004), Nonparametric independent component analysis.
\newblock \emph{Bernoulli}, \textbf{10}, 565-582.

\end{thebibliography}
\end{document}